\documentclass[12 pt]{amsart}
\usepackage{amssymb,float}
\usepackage[colorlinks=true]{hyperref}
\usepackage[all]{xy}
\usepackage[toc,page,title,titletoc,header]{appendix}
\usepackage{xcolor}
\usepackage{tikz-cd}
\usepackage{graphicx}
\usepackage{epigraph}
\usepackage{subcaption}

\usepackage{svg}

\begin{document}
	\pdfoutput=1
	\theoremstyle{plain}
	\newtheorem{thm}{Theorem}[section]
	\newtheorem*{thm1}{Theorem 1}
	
	\newtheorem*{thmM}{Main Theorem}
\newtheorem{thmA}{Theorem}
\renewcommand{\thethmA}{A}
    \theoremstyle{plain}
    \newtheorem{thmC}{Theorem}
    \renewcommand{\thethmC}{C}
	\newtheorem*{thm2}{Theorem 2}
	\newtheorem{lemma}[thm]{Lemma}
	\newtheorem{lem}[thm]{Lemma}
	\newtheorem{cor}[thm]{Corollary}
	\newtheorem{pro}[thm]{Proposition}
	\newtheorem{prop}[thm]{Proposition}
	\newtheorem{variant}[thm]{Variant}
	\theoremstyle{definition}
	\newtheorem{notations}[thm]{Notations}
	\newtheorem{rem}[thm]{Remark}
	\newtheorem{rmk}[thm]{Remark}
	\newtheorem{rmks}[thm]{Remarks}
	\newtheorem{defi}[thm]{Definition}
	\newtheorem{exe}[thm]{Example}
	\newtheorem{claim}[thm]{Claim}
	\newtheorem{ass}[thm]{Assumption}
	\newtheorem{prodefi}[thm]{Proposition-Definition}
	\newtheorem{que}[thm]{Question}
	\newtheorem{con}[thm]{Conjecture}
	
	\newtheorem{exa}[thm]{Example}
	\newtheorem*{assa}{Assumption A}

	\newtheorem*{dmlcon}{Dynamical Mordell-Lang Conjecture}
	\newtheorem*{condml}{Dynamical Mordell-Lang Conjecture}
	\newtheorem*{congb}{Geometric Bogomolov Conjecture}
	\newtheorem*{congdaocurve}{Dynamical Andr\'e-Oort Conjecture for curves}
	
	\newtheorem*{pdd}{P(d)}
	\newtheorem*{bfd}{BF(d)}

	\newtheorem*{probreal}{Realization problems}
	\numberwithin{equation}{section}
	\newcounter{elno}                
	\def\points{\list
		{\hss\llap{\upshape{(\roman{elno})}}}{\usecounter{elno}}}
	\let\endpoints=\endlist
	\newcommand{\SH}{\rm SH}
	\newcommand{\Cov}{\rm Cov}
	\newcommand{\Tan}{\rm Tan}
	\newcommand{\res}{\rm res}
	\newcommand{\Om}{\Omega}
	\newcommand{\om}{\omega}
	\newcommand{\La}{\Lambda}
	\newcommand{\la}{\lambda}
	\newcommand{\mc}{\mathcal}
	\newcommand{\mb}{\mathbb}
	\newcommand{\surj}{\twoheadrightarrow}
	\newcommand{\inj}{\hookrightarrow}
	\newcommand{\zar}{{\rm zar}}
	\newcommand{\Exc}{{\rm Exc}}
	\newcommand{\Mod}{{\rm Mod}}
	\newcommand{\an}{{\rm an}}
	\newcommand{\red}{{\rm red}}
	\newcommand{\codim}{{\rm codim}}
	\newcommand{\Supp}{{\rm Supp\;}}
	\newcommand{\Leb}{{\rm Leb}}
	\newcommand{\rank}{{\rm rank}}
	\newcommand{\Ker}{{\rm Ker \ }}
	\newcommand{\Pic}{{\rm Pic}}
	\newcommand{\Der}{{\rm Der}}
	\newcommand{\Div}{{\rm Div}}
	\newcommand{\Hom}{{\rm Hom}}
	\newcommand{\Corr}{{\rm Corr}}
	\newcommand{\im}{{\rm im}}
	\newcommand{\Spec}{{\rm Spec \,}}
	\newcommand{\Nef}{{\rm Nef \,}}
	\newcommand{\Frac}{{\rm Frac \,}}
	\newcommand{\Sing}{{\rm Sing}}
	\newcommand{\sing}{{\rm sing}}
	\newcommand{\reg}{{\rm reg}}
	\newcommand{\Char}{{\rm char\,}}
	\newcommand{\Tr}{{\rm Tr}}
	\newcommand{\ord}{{\rm ord}}
	\newcommand{\bif}{{\rm bif}}
	\newcommand{\AS}{{\rm AS}}
	\newcommand{\FS}{{\rm FS}}
	\newcommand{\CE}{{\rm CE}}
	\newcommand{\PCE}{{\rm PCE}}
	\newcommand{\WR}{{\rm WR}}
	\newcommand{\PR}{{\rm PR}}
	\newcommand{\TCE}{{\rm TCE}}
	\newcommand{\diam}{{\rm diam\,}}
	\newcommand{\id}{{\rm id}}
	\newcommand{\NE}{{\rm NE}}
	\newcommand{\Gal}{{\rm Gal}}
	\newcommand{\Min}{{\rm Min \ }}
	\newcommand{\Hol}{{\rm Hol \ }}
	\newcommand{\Rat}{{\rm Rat}}
	\newcommand{\FL}{{\rm FL}}
	\newcommand{\fm}{{\rm fm}}
	
	\newcommand{\Max}{{\rm Max \ }}
	\newcommand{\Alb}{{\rm Alb}\,}
	\newcommand{\Aff}{{\rm Aff}\,}
	\newcommand{\GL}{{\rm GL}\,}        
	\newcommand{\PGL}{{\rm PGL}\,}
	\newcommand{\Bir}{{\rm Bir}}
	\newcommand{\Aut}{{\rm Aut}}
	\newcommand{\End}{{\rm End}}
	\newcommand{\Bel}{{\rm Bel}}
	\newcommand{\QC}{{\rm QC}}

	\newcommand{\Per}{{\rm Per}\,}
	\newcommand{\Preper}{{\rm Preper}\,}
	\newcommand{\ie}{{\it i.e.\/},\ }
	\newcommand{\niso}{\not\cong}
	\newcommand{\nin}{\not\in}
	\newcommand{\soplus}[1]{\stackrel{#1}{\oplus}}
	\newcommand{\by}[1]{\stackrel{#1}{\rightarrow}}
	\newcommand{\longby}[1]{\stackrel{#1}{\longrightarrow}}
	\newcommand{\vlongby}[1]{\stackrel{#1}{\mbox{\large{$\longrightarrow$}}}}
	\newcommand{\ldownarrow}{\mbox{\Large{\Large{$\downarrow$}}}}
	\newcommand{\lsearrow}{\mbox{\Large{$\searrow$}}}
	\renewcommand{\d}{\stackrel{\mbox{\scriptsize{$\bullet$}}}{}}
	\newcommand{\dlog}{{\rm dlog}\,}    
	\newcommand{\longto}{\longrightarrow}
	\newcommand{\vlongto}{\mbox{{\Large{$\longto$}}}}
	\newcommand{\limdir}[1]{{\displaystyle{\mathop{\rm lim}_{\buildrel\longrightarrow\over{#1}}}}\,}
	\newcommand{\liminv}[1]{{\displaystyle{\mathop{\rm lim}_{\buildrel\longleftarrow\over{#1}}}}\,}
	\newcommand{\norm}[1]{\mbox{$\parallel{#1}\parallel$}}
	\newcommand{\boxtensor}{{\Box\kern-9.03pt\raise1.42pt\hbox{$\times$}}}
	\newcommand{\into}{\hookrightarrow}
	\newcommand{\image}{{\rm image}\,}
	\newcommand{\Lie}{{\rm Lie}\,}      
	\newcommand{\CM}{\rm CM}
	\newcommand{\Teich}{\rm Teich\;}
	\newcommand{\sext}{\mbox{${\mathcal E}xt\,$}}  
	\newcommand{\shom}{\mbox{${\mathcal H}om\,$}}  
	\newcommand{\coker}{{\rm coker}\,}  
	\newcommand{\sm}{{\rm sm}}
	\newcommand{\pgcd}{\text{pgcd}}
	\newcommand{\trd}{\text{tr.d.}}
	\newcommand{\tensor}{\otimes}
	\newcommand{\hotimes}{\hat{\otimes}}
	
	\newcommand{\CH}{{\rm CH}}
	\newcommand{\tr}{{\rm tr}}
	\newcommand{\e}{\rm SH}
	
	\renewcommand{\iff}{\mbox{ $\Longleftrightarrow$ }}
	\newcommand{\supp}{{\rm supp}\,}
	\newcommand{\ext}[1]{\stackrel{#1}{\wedge}}
	\newcommand{\onto}{\mbox{$\,\>>>\hspace{-.5cm}\to\hspace{.15cm}$}}
	\newcommand{\propsubset}
	{\mbox{$\textstyle{
				\subseteq_{\kern-5pt\raise-1pt\hbox{\mbox{\tiny{$/$}}}}}$}}
	
	\newcommand{\sA}{{\mathcal A}}
	\newcommand{\sB}{{\mathcal B}}
	\newcommand{\sC}{{\mathcal C}}
	\newcommand{\sD}{{\mathcal D}}
	\newcommand{\sE}{{\mathcal E}}
	\newcommand{\sF}{{\mathcal F}}
	\newcommand{\sG}{{\mathcal G}}
	\newcommand{\sH}{{\mathcal H}}
	\newcommand{\sI}{{\mathcal I}}
	\newcommand{\sJ}{{\mathcal J}}
	\newcommand{\sK}{{\mathcal K}}
	\newcommand{\sL}{{\mathcal L}}
	\newcommand{\sM}{{\mathcal M}}
	\newcommand{\sN}{{\mathcal N}}
	\newcommand{\sO}{{\mathcal O}}
	\newcommand{\sP}{{\mathcal P}}
	\newcommand{\sQ}{{\mathcal Q}}
	\newcommand{\sR}{{\mathcal R}}
	\newcommand{\sS}{{\mathcal S}}
	\newcommand{\sT}{{\mathcal T}}
	\newcommand{\sU}{{\mathcal U}}
	\newcommand{\sV}{{\mathcal V}}
	\newcommand{\sW}{{\mathcal W}}
	\newcommand{\sX}{{\mathcal X}}
	\newcommand{\sY}{{\mathcal Y}}
	\newcommand{\sZ}{{\mathcal Z}}
	
	\newcommand{\A}{{\mathbb A}}
	\newcommand{\B}{{\mathbb B}}
	\newcommand{\C}{{\mathbb C}}
	\newcommand{\D}{{\mathbb D}}
	\newcommand{\E}{{\mathbb E}}
	\newcommand{\F}{{\mathbb F}}
	\newcommand{\G}{{\mathbb G}}
	\newcommand{\HH}{{\mathbb H}}
	\newcommand{\LL}{{\mathbb L}}
	\newcommand{\J}{{\mathbb J}}
	\newcommand{\M}{{\mathbb M}}
	\newcommand{\N}{{\mathbb N}}
	\renewcommand{\P}{{\mathbb P}}
	\newcommand{\Q}{{\mathbb Q}}
	\newcommand{\R}{{\mathbb R}}
	\newcommand{\T}{{\mathbb T}}
	\newcommand{\U}{{\mathbb U}}
	\newcommand{\V}{{\mathbb V}}
	\newcommand{\W}{{\mathbb W}}
	\newcommand{\X}{{\mathbb X}}
	\newcommand{\Y}{{\mathbb Y}}
	\newcommand{\Z}{{\mathbb Z}}
	
	\newcommand{\bk}{{\mathbf{k}}}
	
	\newcommand{\bp}{{\mathbf{p}}}
	\newcommand{\ep}{\varepsilon}
	\newcommand{\bbk}{{\overline{\mathbf{k}}}}
	\newcommand{\Fix}{\mathrm{Fix}}
	
	\newcommand{\tor}{{\mathrm{tor}}}
	\renewcommand{\div}{{\mathrm{div}}}
	
	\newcommand{\trdeg}{{\mathrm{trdeg}}}
	\newcommand{\Stab}{{\mathrm{Stab}}}
	
	\newcommand{\OK}{{\overline{K}}}
	\newcommand{\ok}{{\overline{k}}}
	
	\newcommand{\cf}{{\color{red} [c.f. ?]}}
	\newcommand{\jy}{\color{red} jy:}

	\newcommand{\Chat}{\widehat{\C}}
	\newcommand{\Rhat}{\widehat{\R}}
	\newcommand{\Mob}{{\rm Mob}\,}
	\newcommand{\AZ}{Avanzi--Zannier}
	\newcommand{\Muller}{M\"uller}
	\newcommand{\eps}{\varepsilon}
	\newcommand{\crit}{{\rm Crit}\,}
	\newcommand{\orb}{{\rm Orb}\,}
	\newcommand{\mult}{{\rm mult}\,}
	\newcommand{\ordd}{{\rm ord}\,}
	\theoremstyle{plain}
    \newtheorem{thmB}{Theorem}
    \renewcommand{\thethmB}{B}

	\title[]{Invariant real curves arising from semiconjugacies of rational maps}

	\author{Hongxin Sun}

\address{School of Mathematical Sciences, Fudan University, Shanghai 200433, China}

\address{Institute for Theoretical Sciences, Westlake University, Hangzhou 310030, China}

	\email{sunhongxin@westlake.edu.cn}
	
	\date{\today}

	\bibliographystyle{amsplain}
	
	
\begin{abstract}
In this paper, we study real analytic invariant curves arising from semiconjugacies of real rational maps. Such semiconjugacies arise naturally in Eremenko's work on invariant Jordan analytic curves. While invariant generalized circles are easy to obtain from real rational maps, we are interested in curves not contained in any generalized circle. First, for a non-exceptional polynomial of degree at least two, we prove that every invariant curve arising from such a semiconjugacy is contained in a generalized circle. Second, in the rational map case, we construct invariant curves arising in this way that are not contained in any generalized circle and have arbitrarily many self-intersections, counted with multiplicity. Moreover, the rational maps can be chosen so that their restrictions to these curves have positive topological entropy.
\end{abstract}

	\maketitle

	\section{Introduction}
    \subsection{Invariant curves and semiconjugacies}

The study of analytic invariant curves for rational maps goes back to Fatou. In the final paragraph of his memoir on functional equations, Fatou noted that analytic curves invariant under rational transformations remained to be studied \cite[p.~314]{Fatou1920}.

If \(f\) is a real rational map, then \(\Rhat\) is \(f\)-invariant. A generalized circle in \(\Chat\) is a M\"obius image of \(\Rhat\). Rational maps preserving a given generalized circle are precisely the corresponding M\"obius conjugates of real rational maps.

We are interested in invariant curves not contained in any generalized circle. Rotation domains provide another source of invariant curves. On a Siegel disk or a Herman ring, a suitable iterate of \(f\) is conformally conjugate to an irrational rotation of a disk or an annulus, respectively. The inverse images of level circles of positive radius in the rotation domain are real analytic curves invariant under that iterate. However, these curves contain no periodic points of \(f\), in contrast to the setting considered below \cite{EremenkoInvariant}.

Suppose that an invariant real analytic curve \(\Gamma\) is smooth near a repelling periodic point \(a\) of period \(q\). The derivative of \(f^q\) preserves the tangent line to \(\Gamma\) at \(a\), so its multiplier \(\lambda=(f^q)'(a)\) is real. This observation is used in Eremenko's study of invariant curves and semiconjugacies \cite{EremenkoInvariant}.

In particular, Eremenko proved that if \(\Gamma\subset J(f)\) is an invariant Jordan analytic curve containing a repelling periodic point and \(f|_\Gamma\) is not a homeomorphism, then \(\Gamma\) is a real semi-algebraic curve unless \(f\) is a Latt\`es map \cite[Theorem~1 and its corollary]{EremenkoInvariant}. Actually, in the non-Latt\`es case, after replacing \(f\) by an iterate, \(\Gamma\) is obtained from a rational semiconjugacy in the sense of Definition~\ref{def:semiconjugacy} below.

In a related direction, Eremenko and van Strien proved that if a non-Latt\`es rational map has real multipliers at all repelling periodic points, then its Julia set is contained in a generalized circle \cite[Theorem~1]{EVS}.

In this paper, we study invariant curves arising from rational semiconjugacies.

\begin{defi}
\label{def:semiconjugacy}
Let \(f\in\C(z)\) be a nonconstant rational map. A real analytic curve \(\Gamma\subset\Chat\) is called an \emph{invariant curve of \(f\) arising from a rational semiconjugacy} if there exist nonconstant maps \(g\in\R(z)\) and \(h\in\C(z)\) such that
\[
    h\circ g=f\circ h,
    \qquad
    \Gamma=h(\Rhat).
\]
We refer to \((g,f,h,\Gamma)\) as the corresponding semiconjugacy quadruple.

We have \(f(\Gamma)\subset\Gamma\), and the following diagram commutes:
\[
\begin{tikzcd}
    \Rhat \arrow[r,"g"] \arrow[d,"h"'] &
    \Rhat \arrow[d,"h"] \\
    \Gamma \arrow[r,"f"'] &
    \Gamma .
\end{tikzcd}
\]
\end{defi}

\subsection{Polynomial case}

We use the normalized Chebyshev polynomials, characterized by \( T_n(z+z^{-1})=z^n+z^{-n}. \) A rational map of degree \(d\ge2\) is called \emph{exceptional} if it is a Latt\`es map or is M\"obius conjugate to \(z^{\pm d}\) or \(\pm T_d\); see \cite[Lemma~2]{EVS}. Otherwise, it is called \emph{non-exceptional}.

Our first result concerns non-exceptional polynomials \(f\), while allowing \(g\) and \(h\) to be rational maps.

\begin{thmA}\label{thm:A-nonexceptional}
Let \(f\in\C[z]\) be a non-exceptional polynomial of degree at least two. Then every invariant curve of \(f\) arising from a rational semiconjugacy in the sense of Definition~\ref{def:semiconjugacy} is contained in a generalized circle.
\end{thmA}

The non-exceptionality assumption cannot be omitted, as the following example with a Chebyshev polynomial shows. Let \(\alpha=e^{2\pi i/3}\), and define
\[
    f=T_4,\qquad g(z)=z^4,\qquad
    h(z)=\alpha z+\frac{1}{\alpha z}.
\]
Since \(\alpha^3=1\), the identity \(T_4(w+w^{-1})=w^4+w^{-4}\) gives \(h\circ g=f\circ h\). Thus \(\Gamma:=h(\Rhat)\) is an invariant curve of \(f\) arising from a rational semiconjugacy.

However, \(\Gamma\) is a hyperbola, which is not contained in any generalized circle.

\subsection{Positive entropy}
\label{subsec:main-entropy}

We first consider special semiconjugacies arising from functional decompositions, as in \cite{EremenkoInvariant,Muller}. Let \(u,v\in\C(z)\) be nonconstant rational maps such that \(g:=v\circ u\in\R(z)\), and set
\[
    f:=u\circ v,
    \qquad
    \Gamma:=u(\Rhat).
\]
Thus \(u\circ g=f\circ u\). Since \(g\) preserves \(\Rhat\), we obtain the following commutative diagram by restriction:
\[
\begin{tikzcd}
    \Rhat \arrow[r,"g=v\circ u"] \arrow[d,"u"'] &
    \Rhat \arrow[d,"u"] \\
    \Gamma \arrow[r,"f=u\circ v"'] &
    \Gamma .
\end{tikzcd}
\]
The following theorem gives a condition for obtaining positive entropy without changing \(\Gamma\).

\begin{thmB}\label{thm:entropy-main}
Let \(u,v\in\C(z)\) be nonconstant rational maps such that \(g:=v\circ u\in\R(z)\), and set \(\Gamma:=u(\Rhat)\). Suppose that there exists \(a\in\Rhat\) such that
\[
    n:=\#\{x\in\Rhat:g(x)=a\}\ge2,
\]
and every point of this real fiber is non-critical for \(g\). Then there exists a real M\"obius transformation \(A\) such that, with
\[
    g_A:=A\circ g,
    \qquad
    f_A:=u\circ A\circ v,
\]
we have
\[
    u\circ g_A=f_A\circ u,
    \qquad
    f_A(\Gamma)\subset\Gamma,
\]
and
\[
    h_{\mathrm{top}}(f_A|_\Gamma)
    =h_{\mathrm{top}}(g_A|_{\Rhat})
    \ge\log n>0.
\]
\end{thmB}

\subsection{Rational map case}
\label{subsec:main-rational-examples}

For non-exceptional polynomials, Theorem~\ref{thm:A-nonexceptional} shows that invariant curves arising from rational semiconjugacies are contained in generalized circles. For rational maps, the following result gives a different conclusion.

\begin{thmC}\label{thm:rational-main}
For every integer \(N\ge1\), there exists a semiconjugacy quadruple \((g,f,h,\Gamma)\) as in Definition~\ref{def:semiconjugacy}, with \(\deg f=\deg g\ge2\), such that \(\Gamma\) is not contained in any generalized circle and has at least \(N\) self-intersections, counted with multiplicity. Moreover, \(h_{\mathrm{top}}(f|_\Gamma)>0\).
\end{thmC}

Self-intersections are counted on \(\Chat\), including at infinity. An ordinary \(k\)-fold point contributes \(\binom{k}{2}\) to the self-intersection count. In fact, our construction produces a self-intersection point of arbitrarily large multiplicity, rather than arbitrarily many distinct self-intersection points.

\subsection{Strategy of the proof}

For Theorem~\ref{thm:A-nonexceptional}, we first use non-ex\-cep\-tion\-al\-ity to reduce rational semiconjugacies to polynomial ones without changing the curve. We then apply the polynomial decomposition theory of Ritt and Engstrom \cite{Ritt,Engstrom}, in Inou's formulation for semiconjugacies \cite{Inou}, together with the real decomposition lemma stated by M\"uller \cite[Lemma~5.1]{Muller}. These results allow us to argue by induction on the degree of the semiconjugating map. We also use B\"ottcher coordinates to finish the proof.

For Theorem~\ref{thm:entropy-main}, we change the dynamics without changing the curve. We choose a real M\"obius transformation \(A\) so that \(G:=A\circ g\) is a \emph{Bernoulli map}. Here this means that a closed interval \(I\subset\Rhat\) contains \(n\ge2\) pairwise disjoint closed subintervals \(I_1,\ldots,I_n\), each mapped bijectively onto \(I\) by \(G\). Such maps have positive entropy on \(\Rhat\); see Proposition~\ref{prop:full-branch}. It follows from the finite-to-one semiconjugacy that the entropy on \(\Gamma\) is positive.

The examples for Theorem~\ref{thm:rational-main} are obtained from the construction of Avanzi and Zannier \cite[Proposition~4.7(3)]{AZ}, which is also discussed by M\"uller \cite[Example~5.3]{Muller}. In this construction, arbitrarily many branches can meet at infinity. We prove that their tangent lines are distinct, so that the resulting self-intersection point has arbitrarily large multiplicity. Theorem~\ref{thm:entropy-main} then yields positive entropy while leaving the curve unchanged.
\subsection{Organization of the paper}

In Section~\ref{sec:prelim}, we introduce the notation and recall some basic facts. Theorems~\ref{thm:A-nonexceptional}, \ref{thm:entropy-main}, and \ref{thm:rational-main} are proved in Sections~\ref{sec:poly}, \ref{sec:entropy}, and \ref{sec:rational}, respectively.

\subsection{Acknowledgements}

The author is grateful to her advisor, Zhuchao Ji, for his guidance, helpful discussions, and suggestions on the writing of the manuscript. The author also thanks Mai Hao for useful discussions and Alexandre Eremenko for helpful discussions.
\subsection{AI disclosure}

All statements and proofs in the paper were developed by the author.  GPT was used to check and edit the paper after a draft was written.

\section{Preliminaries}\label{sec:prelim}

We write \(\Chat:=\C\cup\{\infty\}\) and \(\Rhat:=\R\cup\{\infty\}\). Let \(\Mob(\Chat)\) denote the group of M\"obius transformations of \(\Chat\), and let \(\Aff(\C)\) denote the group of maps \(z\mapsto az+b\), where \(a\in\C^*\) and \(b\in\C\). For \(U\in\C(z)\), let \(\overline U\) denote the rational function obtained by conjugating its coefficients, so that \(\overline U(z)=\overline{U(\overline z)}\). For nonconstant \(U\in\C(z)\), we have
\begin{equation}\label{eq:real-rational}
    U\in\R(z)
    \quad\Longleftrightarrow\quad
    U=\overline U
    \quad\Longleftrightarrow\quad
    U(\Rhat)\subset\Rhat.
\end{equation}

A \emph{generalized circle} in \(\Chat\) is a M\"obius image of \(\Rhat\). It is either a Euclidean circle or a set \(L\cup\{\infty\}\), where \(L\subset\C\) is an affine line.

\begin{defi}[Real analytic curve]
We use a parametrized notion. A set \(\Gamma\subset\Chat\) is called a real analytic curve if there exist an interval \([a,b]\) and a real analytic map \(\phi:[a,b]\to\Chat\) such that \(\Gamma=\phi([a,b])\) and \(\phi\) extends real analytically to \((a-\varepsilon,b+\varepsilon)\) for some \(\varepsilon>0\).
\end{defi}

Note that we do not require \(\phi\) to be injective.

We shall use the following elementary criterion.

\begin{lemma}\label{lem:circle-real}
Let \(U\in\C(z)\) be nonconstant. Then \(U(\Rhat)\) is contained in a generalized circle if and only if there exists \(\lambda\in\Mob(\Chat)\) such that \(\lambda\circ U\in\R(z)\).
\end{lemma}

\begin{proof}
For every \(\lambda\in\Mob(\Chat)\), \eqref{eq:real-rational} gives
\[
    \lambda\circ U\in\R(z)
    \quad\Longleftrightarrow\quad
    U(\Rhat)\subset\lambda^{-1}(\Rhat).
\]
Every generalized circle has the form \(\lambda^{-1}(\Rhat)\), so the result follows.
\end{proof}

\begin{rem}\label{rem:polynomial-circle}
If \(U\in\C[z]\) is nonconstant, then \(U(\infty)=\infty\). Since every generalized circle containing \(\infty\) has the form \(L\cup\{\infty\}\) for an affine line \(L\subset\C\), \(U(\Rhat)\) is contained in a generalized circle if and only if \(U(\R)\) is contained in an affine line. Equivalently, there exists \(A\in\Aff(\C)\) such that \( A\circ U\in\R[z]. \)
\end{rem}

	\section{Polynomial case}\label{sec:poly}

The following theorem gives the conclusion of Theorem~\ref{thm:A-nonexceptional} for polynomial semiconjugacies, without the non-exceptionality assumption. We first use it to prove Theorem~\ref{thm:A-nonexceptional}, and prove Theorem~\ref{thm:poly-main} at the end of this section.

\begin{thm}\label{thm:poly-main}
Let \(g\in\R[z]\) and \(f,h\in\C[z]\) satisfy
\[
    h\circ g=f\circ h,
    \qquad
    \deg f=\deg g\ge2,
    \qquad
    \deg h\ge1.
\]
Then \(\Gamma:=h(\Rhat)\) is contained in a generalized circle.
\end{thm}

\begin{proof}[Proof of Theorem~\ref{thm:A-nonexceptional} assuming Theorem~\ref{thm:poly-main}]

Choose \(g\in\R(z)\) and \(h\in\C(z)\) as in Definition~\ref{def:semiconjugacy}, so that
\[
    h\circ g=f\circ h,
    \qquad \Gamma=h(\Rhat).
\]
Then \(\deg g=\deg f=:d\ge2\). By \cite[Corollary~4.3]{MilnorLattes}, \(g\) is non-exceptional. The aim is to prove that \(\Gamma\) is contained in a generalized circle. We claim that \(h^{-1}(\infty)=\{p\}\) for some \(p\in\Rhat\).

Define \(P=h^{-1}(\infty)\), the nonempty finite set of poles of \(h\) on \(\Chat\). Since \(f^{-1}(\infty)=\{\infty\}\), the semiconjugacy gives \(g^{-1}(P)=P\). As \(g\) has real coefficients, \(g^{-1}(\overline P)=\overline P\) as well. Here \(\overline P\) is the set consisting of the complex conjugates of all points in \(P\). Hence the nonempty finite set \(S=P\cup\overline P\) satisfies \(g^{-1}(S)=S\).

By the standard description of finite completely invariant sets for rational maps \cite{MilnorIntro}, \(S\) has one or two points. If \(\#S=2\), then \(g\) is M\"obius conjugate to \(z^{\pm d}\), contradicting its non-exceptionality. Hence \(S=\{p\}\). Since \(P\) is nonempty and \(S\) is invariant under complex conjugation, we have \(P=\{p\}\) with \(p\in\Rhat\), proving the claim.

Choose a real M\"obius transformation \(M\) with \(M(p)=\infty\), and define
\[
    \widetilde g=M\circ g\circ M^{-1},
    \qquad \widetilde h=h\circ M^{-1}.
\]
The identities \(g^{-1}(p)=\{p\}\) and \(h^{-1}(\infty)=\{p\}\) imply that infinity is the only pole of each of the two transformed maps. Hence \(\widetilde g\in\R[z]\) and \(\widetilde h\in\C[z]\), with
\[
    \widetilde h\circ\widetilde g=f\circ\widetilde h,
    \qquad
    \widetilde h(\Rhat)=h(\Rhat)=\Gamma.
\]
The last equality uses \(M(\Rhat)=\Rhat\). Theorem~\ref{thm:poly-main} now gives the conclusion.
\end{proof}

We now prove Theorem~\ref{thm:poly-main}.

We use the following elementary fact repeatedly. For \(P\in\C[z]\) and an affine line \(L\subset\C\), \(P(\R)\subset L\) if and only if there exists an affine map \(B\) with \(B(L)=\R\) such that \(B\circ P\in\R[z]\).

We shall also use the following real decomposition lemma, in the form stated by M\"uller.

\begin{lemma}[{\cite[Lemma~5.1]{Muller}}]
\label{lem:real-decomposition}
Let \(V\in\C(z)\) and \(W\in\C[z]\) be nonconstant. If \(V\circ W\in\R(z)\), then there exists an affine map \(A\in\Aff(\C)\) such that \( A\circ W\in\R[z]. \)
\end{lemma}

We call a polynomial semiconjugacy \(h\circ g=f\circ h\) \emph{primitive} when \(\deg f\) and \(\deg h\) are coprime. We first prove Theorem~\ref{thm:poly-main} in the primitive case, and then treat the remaining case by induction on \(\deg h\).

\subsection{The primitive case}

We use the following normal forms obtained from Ritt's second theorem \cite{Ritt}.

\begin{thm}[Primitive normal forms {\cite[Theorem~8]{Inou}}]
\label{thm:inou-primitive}
Let \(f,g,h\in\C[z]\) satisfy
\[
    h\circ g=f\circ h,
    \qquad
    a:=\deg f=\deg g\ge2,
    \qquad
    b:=\deg h>1,
\]
and assume that \(\gcd(a,b)=1\). Then there exist affine maps \(\sigma_1,\sigma_2\) and polynomials \(F,G,H\) such that
\begin{equation}\label{eq:affine-normal-form}
\begin{aligned}
    f&=\sigma_1\circ F\circ\sigma_1^{-1},\\
    g&=\sigma_2\circ G\circ\sigma_2^{-1},\\
    h&=\sigma_1\circ H\circ\sigma_2^{-1},
\end{aligned}
\end{equation}
where one of the following forms occurs.

The power type is
\begin{equation}\label{eq:power-type}
    H(z)=z^b,
    \qquad
    G(z)=z^cP(z^b),
    \qquad
    F(z)=z^cP(z)^b,
\end{equation}
where \(P\in\C[z]\setminus\{0\}\), \(c\ge0\) is an integer, and
\[
    a=c+b\deg P,
    \qquad
    \gcd(c,b)=1.
\]
In the Chebyshev case \cite[Theorem~2.4]{PakovichPolynomial}, we have
\begin{equation}\label{eq:cheb-type-poly}
    G=\eta T_a,
    \qquad
    F=\eta^bT_a,
    \qquad
    H=T_b,
    \qquad
    \eta\in\{1,-1\}.
\end{equation}
Here \(T_n\) denotes the Chebyshev polynomial normalized by \( T_n(z+z^{-1})=z^n+z^{-n}. \)
\end{thm}

To prove the primitive case, we study invariant affine lines in B\"ottcher coordinates. For a polynomial \(R\) of degree \(d\ge2\) and leading coefficient \(A_R\), let \(\Phi_R\) be the normalized B\"ottcher coordinate at infinity:
\begin{equation}\label{eq:normalized-bottcher-coordinate}
    \Phi_R(R(z))=A_R\bigl(\Phi_R(z)\bigr)^d,
    \qquad
    \Phi_R(z)=z+O(1)\quad(z\to\infty).
\end{equation}
Existence and uniqueness follow from \cite[Theorem~6.7]{MilnorIntro}. If \(R\in\R[z]\), uniqueness applied to \(z\mapsto\overline{\Phi_R(\overline z)}\) gives
\begin{equation}\label{eq:real-bottcher-coordinate}
    \Phi_R(\overline z)=\overline{\Phi_R(z)}.
\end{equation}

\begin{lemma}\label{lem:centered-invariant-line}
Let \(G\in\C[z]\) have degree \(d\ge2\), and assume that the coefficient of \(z^{d-1}\) in \(G\) is zero. If an affine line \(L\subset\C\) satisfies \(G(L)\subset L\), then \(0\in L\). Moreover, \( \Phi_G(L)=L \) as germs at infinity.
\end{lemma}

\begin{proof}
Comparing the coefficients of \(z^{d-1}\) in the B\"ottcher identity gives
\begin{equation}\label{eq:power-bottcher-no-translation}
    \Phi_G(z)=z+O(z^{-1})
    \qquad(z\to\infty).
\end{equation}
Write \(L=\alpha\R+\beta\), with \(\alpha\ne0\), and define
\[
    s(t):=\alpha t+\beta,
    \qquad
    Q:=s^{-1}\circ G\circ s.
\]
Since \(G(L)\subset L\), we have \(Q\in\R[t]\).

Let \(A_G\) and \(A_Q\) be the leading coefficients of \(G\) and \(Q\), respectively. Then \(A_Q=A_G\alpha^{d-1}\), and
\[
\begin{aligned}
    \Phi_G(G(z))&=A_G\bigl(\Phi_G(z)\bigr)^d,\\
    \Phi_Q(Q(t))&=A_Q\bigl(\Phi_Q(t)\bigr)^d.
\end{aligned}
\]
Define \( \Psi(t):=\alpha^{-1}\Phi_G(s(t)). \) Using \(s\circ Q=G\circ s\), we obtain
\[
\begin{aligned}
    \Psi(Q(t))
    &=\alpha^{-1}\Phi_G(s(Q(t)))\\
    &=\alpha^{-1}\Phi_G(G(s(t)))\\
    &=\alpha^{-1}A_G\bigl(\Phi_G(s(t))\bigr)^d\\
    &=A_G\alpha^{d-1}\Psi(t)^d\\
    &=A_Q\Psi(t)^d.
\end{aligned}
\]
Moreover, \(\Psi\) is locally biholomorphic near infinity and satisfies \(\Psi(t)=t+O(1)\) as \(t\to\infty\). By uniqueness of the normalized B\"ottcher coordinate, \(\Psi=\Phi_Q\). Hence \( \Phi_G\circ s=\alpha\Phi_Q. \)

Since \(Q\in\R[t]\), equation \eqref{eq:real-bottcher-coordinate} gives \(\Phi_Q(t)\in\R\) for real \(t\) of sufficiently large modulus. Thus
\[
    \Phi_G(\alpha t+\beta)=\alpha\Phi_Q(t)\in\alpha\R.
\]
The Euclidean distance between \(L\) and \(\alpha\R\) therefore satisfies
\[
    0\le\operatorname{dist}(L,\alpha\R)
    \le
    \bigl|\Phi_G(\alpha t+\beta)-(\alpha t+\beta)\bigr|.
\]
By \eqref{eq:power-bottcher-no-translation}, the right-hand side tends to zero as \(|t|\to\infty\). Since \(L\) and \(\alpha\R\) are parallel, they coincide. In particular, \(0\in L\).

We have also shown that \(\Phi_G\) maps the germ of \(L\) at infinity into itself. Local invertibility gives \(\Phi_G(L)=L\) as germs at infinity.
\end{proof}

We now apply this lemma to the two forms above.

\begin{lemma}[An invariant-line lemma]\label{lem:line}
Let \(G,H\in\C[z]\), with \(\deg G\ge2\), be of one of the following forms:
\begin{enumerate}
\item
\(H(z)=z^b\) and \(G(z)=z^cP(z^b)\), where \(b>1\) and \(c\ge0\) are integers and \(P\in\C[z]\setminus\{0\}\);
\item
\(G=\eta T_a\) and \(H=T_b\), where \(a\ge2\) and \(b\ge1\) are integers and \(\eta\in\{1,-1\}\).
\end{enumerate}
If an affine line \(L\subset\C\) satisfies \(G(L)\subset L\), then \(H(L)\) is contained in an affine line.
\end{lemma}

\begin{proof}
Let \(d=\deg G\). In both cases, the coefficient of \(z^{d-1}\) in \(G\) is zero. By Lemma~\ref{lem:centered-invariant-line}, \(L=\alpha\R\) for some \(\alpha\ne0\), and \(\Phi_G(L)=L\) as germs at infinity. Recall that \(\Phi_G\) denotes the normalized B\"ottcher coordinate of \(G\) at infinity.

\smallskip
\noindent (1) We immediately obtain \( H(L)\subset\alpha^b\R. \)

\smallskip
\noindent (2) Define
\[
    q(z):=z+z^{-1},
    \qquad
    S(z):=\eta z^a.
\]
Thus
\begin{equation}\label{eq:cheb-power-coordinate}
    G\circ q=q\circ S.
\end{equation}
Since \(\infty\) is not a critical point of \(q\), let \(\psi\) be the local inverse of \(q\) fixing \(\infty\). Since \(q(z)=z+O(z^{-1})\), we have \(\psi(w)=w+O(w^{-1})\) as \(w\to\infty\). Equation~\eqref{eq:cheb-power-coordinate} gives
\[
    \psi(G(w))
    =S(\psi(w))
    =\eta\bigl(\psi(w)\bigr)^a.
\]
Since the leading coefficient of \(G=\eta T_a\) is \(\eta\), uniqueness of the normalized B\"ottcher coordinate gives \(\psi=\Phi_G\).

It follows from \(\Phi_G(L)=L\) that \(q(L)=L\) as germs at infinity. Hence, for all sufficiently large \(z\in L\), both \(z\) and \(q(z)\) belong to \(L\). Since \(L\) passes through the origin, we have
\begin{equation}\label{eq:cheb-line-inversion}
    \frac1z=q(z)-z\in L.
\end{equation}
Taking \(z=\alpha t\) with \(t\in\R\) in \eqref{eq:cheb-line-inversion} gives \(\alpha^{-1}\in\alpha\R\), hence \(\alpha^2\in\R\). Therefore \(L=\R\) or \(L=i\R\). Then
\[
    T_b(\R)\subset\R,
    \qquad
    T_b(i\R)\subset i^b\R.
\]

In either case, \(H(L)\) is contained in an affine line.
\end{proof}

\begin{prop}[The primitive case]\label{prop:poly-primitive}
Let \(g\in\R[z]\) and \(f,h\in\C[z]\) satisfy
\[
    h\circ g=f\circ h,
    \qquad
    a:=\deg f=\deg g\ge2,
    \qquad
    b:=\deg h\ge1.
\]
If \(\gcd(a,b)=1\), then \(h(\R)\) is contained in an affine line.
\end{prop}

\begin{proof}
The case \(b=1\) is immediate. For \(b>1\), apply Theorem~\ref{thm:inou-primitive} and set \(L=\sigma_2^{-1}(\R)\). Since \(g\in\R[z]\), we have \(G(L)\subset L\). Lemma~\ref{lem:line} shows that \(H(L)\) lies in an affine line. The same holds for \(h(\R)=\sigma_1(H(L))\), since \(\sigma_1\) is affine.
\end{proof}

\subsection{The non-primitive case}

We use the following consequence of the polynomial decomposition theory of Ritt and Engstrom \cite{Ritt,Engstrom}.

\begin{thm}[Reduction of polynomial semiconjugacies {\cite[Theorem~7]{Inou}}]
\label{thm:inou-reduction}
Let \(f,g,h\in\C[z]\) satisfy
\[
    h\circ g=f\circ h,
    \qquad
    a:=\deg f=\deg g\ge2,
    \qquad
    b:=\deg h>1.
\]
Define \(d:=\gcd(a,b)\), and assume \(d>1\). Then there exist polynomials \(\alpha,\widehat h,f_1\) such that
\begin{equation}\label{eq:inou-reduction}
    h=\alpha\circ\widehat h,
    \qquad
    f_1\circ\widehat h=\widehat h\circ g,
\end{equation}
with
\[
    \deg\widehat h=\frac{b}{d},
    \qquad
    \deg\alpha=d,
    \qquad
    \deg f_1=a.
\]
\end{thm}

The original semiconjugacy and \eqref{eq:inou-reduction} also give \( f\circ\alpha\circ\widehat h =\alpha\circ f_1\circ\widehat h. \) Since \(\widehat h\) is nonconstant, we have \( f\circ\alpha=\alpha\circ f_1. \)

For the case \(b\mid a\), we use the following common-right-factor consequence of Engstrom's theorem \cite{Engstrom}; see also \cite[Theorem~11]{Inou}.
\begin{thm}\label{thm:engstrom-right}
Let \(P_1,P_2,Q_1,Q_2\in\C[z]\) be nonconstant and satisfy
\[
    P_1\circ Q_1=P_2\circ Q_2.
\]
Then there exist polynomials \(W,\widehat Q_1,\widehat Q_2\) such that
\[
    Q_1=\widehat Q_1\circ W,
    \qquad
    Q_2=\widehat Q_2\circ W,
    \qquad
    \deg W=\gcd(\deg Q_1,\deg Q_2).
\]
\end{thm}

We now treat the non-primitive case.

\begin{prop}\label{prop:poly-divisible}
Let \(g\in\R[z]\) and \(f,h\in\C[z]\) satisfy
\[
    h\circ g=f\circ h,
    \qquad
    a:=\deg f=\deg g\ge2,
    \qquad
    b:=\deg h>1.
\]
Define \(d:=\gcd(a,b)\), and assume \(d=b>1\). Then \(h(\R)\) is contained in an affine line.
\end{prop}

\begin{proof}
We apply Theorem~\ref{thm:engstrom-right} to \( f\circ h=h\circ g. \) There is a common right factor \(W\) of \(h\) and \(g\) with \( \deg W=\gcd(\deg h,\deg g)=b. \) Since \(W\) is a common right factor of \(h\) and \(g\), there exist nonconstant polynomials \(H,V\) such that
\[
    h=H\circ W,
    \qquad
    g=V\circ W.
\]
The equality \(\deg h=\deg W\) gives \(\deg H=1\). Thus \(\ell:=H^{-1}\in\Aff(\C)\), and \(W=\ell\circ h\). Consequently, \(g=V\circ W=V\circ(\ell\circ h)\). Since \(g\in\R[z]\), Lemma~\ref{lem:real-decomposition} gives an affine map \(A\in\Aff(\C)\) such that \( A\circ W=A\circ\ell\circ h\in\R[z]. \) Hence \(h(\R)\) is contained in the affine line \( (A\circ\ell)^{-1}(\R). \)
\end{proof}

\begin{prop}\label{prop:poly-nonprimitive}
Let \(g\in\R[z]\) and \(f,h\in\C[z]\) satisfy
\[
    h\circ g=f\circ h,
    \qquad
    a:=\deg f=\deg g\ge2,
    \qquad
    b:=\deg h>1.
\]
Define \(d:=\gcd(a,b)\), and assume \(1<d<b\). Then \(h(\R)\) is contained in an affine line.
\end{prop}

\begin{proof}
We prove the proposition by induction on \(b=\deg h\). The case \(b=2\) is vacuous. Fix \(b>2\), and assume that the proposition holds for all smaller values of \(b\).

By Theorem~\ref{thm:inou-reduction}, there exist polynomials \(\alpha,\widehat h,f_1\) such that
\[
    h=\alpha\circ\widehat h,
    \qquad
    f_1\circ\widehat h=\widehat h\circ g,
    \qquad
    f\circ\alpha=\alpha\circ f_1,
\]
with
\[
    \deg\widehat h=\frac{b}{d}<b,
    \qquad
    \deg\alpha=d<b,
    \qquad
    \deg f_1=a.
\]

Since \(\deg\widehat h=b/d<b\), the semiconjugacy \(\widehat h\circ g=f_1\circ\widehat h\) is covered by Proposition~\ref{prop:poly-primitive} if \(\gcd(a,b/d)=1\), by Proposition~\ref{prop:poly-divisible} if \(\gcd(a,b/d)=b/d\), and by the induction hypothesis otherwise. Thus \(\widehat h(\R)\) is contained in an affine line \(L\). Choose an affine map \(B\) such that \(B(L)=\R\). Since
\[
    f_1(\widehat h(\R))
    =\widehat h(g(\R))
    \subset\widehat h(\R)\subset L,
\]
the polynomial \( \widetilde f_1:=B\circ f_1\circ B^{-1} \) takes real values on the infinite subset \(B(\widehat h(\R))\subset\R\). Hence \( \widetilde f_1\in\R[z]. \)

Now define \(H_2:=\alpha\circ B^{-1}\). Then \(\deg H_2=d<b\), and
\[
\begin{aligned}
    f\circ H_2
    &=f\circ\alpha\circ B^{-1}\\
    &=\alpha\circ f_1\circ B^{-1}\\
    &=\alpha\circ B^{-1}\circ B\circ f_1\circ B^{-1}\\
    &=H_2\circ\widetilde f_1.
\end{aligned}
\]
Thus we have another semiconjugacy \( f\circ H_2=H_2\circ\widetilde f_1, \) where
\[
    \widetilde f_1\in\R[z],
    \qquad
    \deg\widetilde f_1=\deg f,
    \qquad
    \deg H_2<b.
\]

Since \(\deg H_2=d<b\) and \(\gcd(a,d)=d>1\), Proposition~\ref{prop:poly-divisible} applies to \(H_2\circ\widetilde f_1=f\circ H_2\). It follows that \(H_2(\R)\) is contained in an affine line.

Finally, since \(B(\widehat h(\R))\subset\R\), we have
\[
\begin{aligned}
    h(\R)
    &=\alpha(\widehat h(\R))\\
    &=(\alpha\circ B^{-1})
        \bigl(B(\widehat h(\R))\bigr)\\
    &\subset(\alpha\circ B^{-1})(\R)
    =H_2(\R).
\end{aligned}
\]
Therefore \(h(\R)\) is contained in an affine line. This completes the induction.
\end{proof}

\begin{proof}[Proof of Theorem~\ref{thm:poly-main}]
Define \(a:=\deg f=\deg g\), \(b:=\deg h\), and \(d:=\gcd(a,b)\). If \(d=1\), Proposition~\ref{prop:poly-primitive} applies. If \(d=b>1\), Proposition~\ref{prop:poly-divisible} applies. If \(1<d<b\), Proposition~\ref{prop:poly-nonprimitive} applies. Thus \(h(\R)\) is contained in an affine line.

Since \(h\) is a nonconstant polynomial, \(h(\infty)=\infty\). Therefore \(\Gamma=h(\Rhat)\) is contained in a generalized circle.
\end{proof}

\section{Positive entropy on invariant curves}
\label{sec:entropy}

In this section we prove Theorem~\ref{thm:entropy-main}. We first establish an entropy bound for Bernoulli maps. Under the assumptions of the theorem, we obtain such a map by postcomposing \(g\) with a suitable real M\"obius transformation.
\subsection{Bernoulli maps and entropy}

We recall the definition of topological entropy, see  \cite[Section~1]{Bowen1971}.

\begin{defi}[Topological entropy]\label{def:topological-entropy}
Let \((X,d)\) be a compact metric space, and let \(G:X\to X\) be a continuous self-map. For each integer \(m\ge1\), define \(d_m(x,y):=\max_{0\le k<m}d(G^k(x),G^k(y))\). Given \(\varepsilon>0\), a subset \(E\subset X\) is said to be \emph{\((m,\varepsilon)\)-separated} if \(d_m(x,y)>\varepsilon\) for any two distinct points \(x,y\in E\). Let \(s_m(\varepsilon)\) be the largest cardinality of such a set. The \emph{topological entropy} of \(G\) is defined by
\begin{equation}\label{eq:topological-entropy}
    h_{\mathrm{top}}(G)
    :=\lim_{\varepsilon\to0^+}
      \limsup_{m\to\infty}
      \frac{1}{m}\log s_m(\varepsilon).
\end{equation}

\end{defi}

We first isolate the dynamical condition that produces the entropy lower bound.  Here and below, by a closed interval in \(\Rhat\) we mean a proper closed arc of the circle \(\Rhat\).

\begin{prop}\label{prop:full-branch}

Let \(G:\Rhat\to\Rhat\) be a rational map, and let \(I\subset\Rhat\) be a closed interval. Suppose that there are a finite index set \(\mathcal I\), with \(|\mathcal I|=n\ge1\), and pairwise disjoint closed intervals \(I_i\subset I\), \(i\in\mathcal I\), such that each restriction \( G|_{I_i}:I_i\longrightarrow I \) is bijective. Then \( h_{\mathrm{top}}(G|_{\Rhat})\ge\log n. \)
\end{prop}

\begin{proof}

The assertion is immediate for \(n=1\). Assume now that \(n\ge2\). Choose a metric \(d\) inducing the topology of \(\Rhat\), and define
\[
    \delta:=\min_{\substack{i,j\in\mathcal I\\i\ne j}}
    \operatorname{dist}(I_i,I_j)>0.
\]

For each \(i\in\mathcal I\), the restriction \(G|_{I_i}\) is a homeomorphism. Denote its inverse by \(\phi_i:I\to I_i\). Fix \(y\in I\). For \(m\ge1\) and each word \(\boldsymbol{i}=(i_0,\ldots,i_{m-1})\in\mathcal I^m\), define
\[
    x_{\boldsymbol{i}}
    :=\phi_{i_0}\circ\cdots\circ\phi_{i_{m-1}}(y),
    \qquad
    E_m:=\{x_{\boldsymbol{i}}:\boldsymbol{i}\in\mathcal I^m\}.
\]
These compositions are well defined since \(\phi_i(I)\subset I\). For each \(0\le k<m\), we have \(G^k(x_{\boldsymbol{i}})\in I_{i_k}\). If \(\boldsymbol{i}\ne\boldsymbol{j}\), then \(i_k\ne j_k\) for some \(0\le k<m\). For this \(k\), \( d\bigl(G^k(x_{\boldsymbol{i}}),G^k(x_{\boldsymbol{j}})\bigr) \ge\delta. \) Thus distinct words give distinct points, and \(E_m\) is an \((m,\delta/2)\)-separated set of cardinality \(n^m\). It follows from \eqref{eq:topological-entropy} that
\[
    h_{\mathrm{top}}(G|_{\Rhat})
    \ge\limsup_{m\to\infty}\frac1m\log s_m(\delta/2)
    \ge\log n.
\]
\end{proof}

\begin{prop}\label{prop:regular-real-preimages}

Let \(G:\Rhat\to\Rhat\) be a nonconstant rational map. Suppose that there exists \(a\in\Rhat\) such that
\[
    n:=\#\{x\in\Rhat:G(x)=a\}\ge2
\]
and none of these \(n\) points is a critical point of \(G\). Then there exists a real M\"obius transformation \(A:\Rhat\to\Rhat\) for which \(h_{\mathrm{top}}\bigl((A\circ G)|_{\Rhat}\bigr)\ge\log n\).
\end{prop}

\begin{proof}
All inverse images in this proof are taken in \(\Rhat\). Write \(G^{-1}(a)=\{x_1,\ldots,x_n\}\), and choose \(p\in\Rhat\setminus G^{-1}(a)\). By the inverse function theorem in local coordinates, there are open intervals \(U_i\ni x_i\), with pairwise disjoint closures avoiding \(p\), such that each \(G|_{U_i}\) is a homeomorphism onto a neighborhood of \(a\). The set \(K_0:=\Rhat\setminus\bigcup_i U_i\) is compact, and \(a\notin G(K_0)\). We may therefore choose a nondegenerate closed interval \(\widetilde I\) such that
\[
    a\in\operatorname{int}(\widetilde I),
    \qquad
    \widetilde I\subset
    \left(\bigcap_iG(U_i)\right)\setminus G(K_0).
\]
Define \(I_i:=(G|_{U_i})^{-1}(\widetilde I)\). Then \( G^{-1}(\widetilde I)=\bigsqcup_{i=1}^n I_i, \) and \(G\) maps each \(I_i\) homeomorphically onto \(\widetilde I\).

Since \(p\notin\bigcup_i I_i\), choose a closed interval \(I\subset\Rhat\setminus\{p\}\) with \(\bigcup_i I_i\subset\operatorname{int}(I)\). Let \(\alpha_-,\alpha_+\) and \(\beta_-,\beta_+\) be the endpoints of \(\widetilde I\) and \(I\), respectively, and choose \(b\in\operatorname{int}(I)\). There is a real M\"obius transformation \(A\) with
\[
    A(\alpha_-)=\beta_-,
    \qquad A(a)=b,
    \qquad A(\alpha_+)=\beta_+.
\]
It maps the open arc between \(\alpha_-\) and \(\alpha_+\) containing \(a\) onto the corresponding arc containing \(b\). Hence \(A(\widetilde I)=I\).

Define \(G_A:=A\circ G\). Then
\[
    G_A^{-1}(I)=G^{-1}(\widetilde I)
    =\bigsqcup_{i=1}^n I_i,
\]
and each \(G_A|_{I_i}:I_i\to I\) is bijective. By Proposition~\ref{prop:full-branch}, \(h_{\mathrm{top}}(G_A|_{\Rhat})\ge\log n\).

\end{proof}

We now incorporate this postcomposition into the semiconjugacy without changing \(\Gamma\) or \(u\).

\begin{prop}\label{prop:modified-semiconjugacy}
Let \(u,v\in\C(z)\) be nonconstant rational maps such that \(g:=v\circ u\in\R(z)\), and define \(\Gamma:=u(\Rhat)\). For every real M\"obius transformation \(A\), define
\[
    g_A:=A\circ v\circ u,
    \qquad
    f_A:=u\circ A\circ v.
\]
Then
\[
    g_A(\Rhat)\subset\Rhat,
    \qquad
    f_A(\Gamma)\subset\Gamma,
\]
and the following diagram commutes:
\[
\begin{tikzcd}
    \Rhat \arrow[r,"g_A"] \arrow[d,"u"'] &
    \Rhat \arrow[d,"u"] \\
    \Gamma \arrow[r,"f_A"'] &
    \Gamma .
\end{tikzcd}
\]
Equivalently, \(u\circ g_A=f_A\circ u\).
\end{prop}

\begin{proof}
Since \(g\in\R(z)\) and \(A\) is real, \(g_A\) preserves \(\Rhat\). We have
\[
    f_A\circ u=u\circ A\circ v\circ u=u\circ g_A.
\]
Consequently,
\[
    f_A(\Gamma)=u\bigl(g_A(\Rhat)\bigr)
    \subset u(\Rhat)=\Gamma.
\]

\end{proof}

The following proposition follows from Bowen's result \cite[Theorem~17]{Bowen1971}, since \(u|_{\Rhat}:\Rhat\to\Gamma\) is a continuous surjection with finite fibers.

\begin{prop}\label{prop:finite-semiconjugacy-entropy}
Let \(u\in\C(z)\) be a nonconstant rational map, and define \(\Gamma:=u(\Rhat)\). Let \(g:\Rhat\to\Rhat\) and \(f:\Gamma\to\Gamma\) be continuous maps satisfying
\[
    u\circ g=f\circ u
    \qquad\text{on }\Rhat.
\]
Then \( h_{\mathrm{top}}(f)=h_{\mathrm{top}}(g). \)
\end{prop}
\subsection{Proof of Theorem~\ref{thm:entropy-main}}

\begin{proof}[Proof of Theorem~\ref{thm:entropy-main}]
By Proposition~\ref{prop:regular-real-preimages}, there is a real M\"obius transformation \(A\) such that
\[
    h_{\mathrm{top}}\bigl((A\circ g)|_{\Rhat}\bigr)\ge\log n.
\]
Define \(g_A:=A\circ g\) and \(f_A:=u\circ A\circ v\). Proposition~\ref{prop:modified-semiconjugacy} gives
\[
    u\circ g_A=f_A\circ u,
    \qquad
    f_A(\Gamma)\subset\Gamma.
\]
Since \(u\) is a nonconstant rational map, Proposition~\ref{prop:finite-semiconjugacy-entropy} yields
\[
    h_{\mathrm{top}}(f_A|_\Gamma)
    =h_{\mathrm{top}}(g_A|_{\Rhat})
    \ge\log n>0.
\]
\end{proof}

\section{Rational map case}\label{sec:rational}
We recall a special case of the construction of Avanzi and Zannier \cite[Proposition~4.7(3)]{AZ}, with \(\ell=1\) in their notation. We use this construction to prove Theorem~\ref{thm:rational-main}. For the self-intersection bound, it suffices to study the branches at infinity. An ordinary \(k\)-fold point has \(k\) smooth local branches with pairwise distinct tangent lines and contributes \(\binom{k}{2}\) to the self-intersection count.

For an integer \(m\ge3\), set
\begin{equation}\label{eq:p1-existence-model}
    v(z)=z(z+1)^m,
    \qquad
    g_1(t)=\frac{1-t}{t^{m+1}-1}
          =-\frac{1}{1+t+\cdots+t^m}.
\end{equation}
Let \(\mu(z)=(z+i)/(z-i)\), and set
\[
    u:=g_1\circ\mu,
    \qquad
    g_0:=v\circ u,
    \qquad
    \Gamma:=u(\Rhat).
\]
Since \(\mu(\Rhat)=S^1\), where \(S^1=\{t\in\C:|t|=1\}\), we have \(\Gamma=g_1(S^1)\).

\begin{prop}\label{prop:p1-existence-model}
For every \(m\ge3\), the map \(g_0\) belongs to \(\R(z)\) and has degree \(m(m+1)\). The set \(\Gamma\) is a real analytic curve with an ordinary \(m\)-fold point at \(\infty\). In particular, it is not contained in any generalized circle and has at least \(\binom{m}{2}\) self-intersections, counted with multiplicity.
\end{prop}

\begin{proof}
Define \(h_1(t)=t^m g_1(t)\). Direct calculation gives
\[
    g_1(1/t)=h_1(t),
    \qquad
    g_1(t)+1=t\bigl(h_1(t)+1\bigr).
\]
Consequently,
\[
    v(g_1(t))
    =g_1(t)t^m\bigl(h_1(t)+1\bigr)^m
    =v(h_1(t)).
\]
Since \(g_1\) has real coefficients and \(\overline\mu=1/\mu\), we obtain \(\overline u=h_1\circ\mu\). As \(v\in\R[z]\), it follows that
\[
    \overline{g_0}
    =v\circ\overline u
    =v\circ h_1\circ\mu
    =v\circ g_1\circ\mu
    =g_0.
\]
Thus \(g_0\in\R(z)\). We have \(\deg g_1=\deg u=m\), whence \(\deg g_0=m(m+1)\). The map \(\theta\mapsto g_1(e^{i\theta})\), viewed as a map into \(\Chat\), is real analytic on \(\R\) and parametrizes \(\Gamma\) on \([0,2\pi]\).

The fiber \(g_1^{-1}(\infty)\) consists of the \(m\) points
\[
    \rho_j=e^{2\pi i j/(m+1)},
    \qquad j=1,\ldots,m.
\]
In the coordinate \(w\mapsto 1/w\) at infinity, the parametrization becomes
\[
    s(t):=\frac{1}{g_1(t)}
         =\frac{t^{m+1}-1}{1-t}
         =-(1+t+\cdots+t^m).
\]
In this coordinate, the tangent vector at \(s=0\) corresponding to \(\rho_j\) is
\[
\begin{aligned}
    \left.\frac{d}{d\theta}s(e^{i\theta})
    \right|_{\theta=2\pi j/(m+1)}
    &=\frac{i(m+1)}{1-\rho_j}\\
    &=-\frac{m+1}{2\sin(\pi j/(m+1))}
       e^{-\pi i j/(m+1)}.
\end{aligned}
\]
These vectors are nonzero and span the pairwise distinct tangent lines
\[
    e^{-\pi i j/(m+1)}\R,
    \qquad j=1,\ldots,m.
\]
Each \(\rho_j\) thus determines a smooth local branch at \(\infty\). Since \(S^1\) is compact and the points \(\rho_j\) are the only preimages of \(\infty\) on \(S^1\), these are all the local branches at \(\infty\). Hence \(\infty\) is an ordinary \(m\)-fold point and contributes \(\binom{m}{2}\) pairs of branches meeting transversely. In particular, \(\Gamma\) is not contained in any generalized circle.
\end{proof}

\begin{proof}[Proof of Theorem~\ref{thm:rational-main}]
Fix \(N\ge1\) and choose \(m\ge3\) such that \(\binom{m}{2}\ge N\). Let \(u,v,g_0,\Gamma\) be as above. By Proposition~\ref{prop:p1-existence-model}, the curve \(\Gamma\) has the required geometric properties. It remains to check the hypothesis of Theorem~\ref{thm:entropy-main}.

The degree \(D:=\deg g_0=m(m+1)\) is even. Since \(g_0(\Rhat)\) is infinite and \(g_0\) has only finitely many critical values, we may choose a finite real value \(a_0\in g_0(\Rhat)\) distinct from all the critical values of \(g_0\). Define
\[
    n:=\#\bigl(g_0^{-1}(a_0)\cap\Rhat\bigr).
\]
The fiber over \(a_0\) in \(\Chat\) consists of \(D\) distinct points, none of which is critical for \(g_0\). Since \(g_0\in\R(z)\), the non-real points in this fiber occur in conjugate pairs. Hence \(n\equiv D\pmod2\). As \(a_0\in g_0(\Rhat)\), we have \(n\ge1\). Thus \(n\ge2\), and none of the points in the real fiber is critical.

We now apply Theorem~\ref{thm:entropy-main} to \(u,v\), with \(g_0=v\circ u\) in place of \(g\). The theorem gives a real M\"obius transformation \(A\) such that the maps
\[
    g:=A\circ g_0=A\circ v\circ u,
    \qquad
    f:=u\circ A\circ v,
\]
satisfy
\[
    u\circ g=f\circ u,
    \qquad
    f(\Gamma)\subset\Gamma,
\]
and
\[
    h_{\mathrm{top}}(f|_\Gamma)
    =h_{\mathrm{top}}(g|_{\Rhat})
    \ge\log n \ge\log 2>0.
\]

Define \(h:=u\). Then \(h\circ g=f\circ h\) and \(\Gamma=h(\Rhat)\). Moreover, \(g\in\R(z)\), \(f,h\in\C(z)\), and \(\deg f=\deg g=D\ge2\). Since \(u\) is unchanged, \(\Gamma\) has the geometric properties established above. Hence \((g,f,h,\Gamma)\) is a quadruple in the sense of Definition~\ref{def:semiconjugacy} with all the required properties.
\end{proof}

\bibliography{bib_reviewed}

@article{AZ,
  author  = {Avanzi, Roberto M. and Zannier, Umberto M.},
  title   = {The equation {$f(X)=f(Y)$} in rational functions
             {$X=X(t), Y=Y(t)$}},
  journal = {Compositio Mathematica},
  volume  = {139},
  number  = {3},
  pages   = {263--295},
  year    = {2003},
  doi     = {10.1023/B:COMP.0000018136.23898.65}
}

@article{Engstrom,
  author  = {Engstrom, H. T.},
  title   = {Polynomial substitutions},
  journal = {American Journal of Mathematics},
  volume  = {63},
  number  = {2},
  pages   = {249--255},
  year    = {1941},
  doi     = {10.2307/2371520}
}

@article{EremenkoInvariant,
  author  = {Eremenko, Alexandre},
  title   = {Invariant curves and semiconjugacies of rational functions},
  journal = {Fundamenta Mathematicae},
  volume  = {219},
  number  = {3},
  pages   = {263--270},
  year    = {2012},
  doi     = {10.4064/fm219-3-5}
}

@article{EVS,
  author  = {Eremenko, Alexandre and van Strien, Sebastian},
  title   = {Rational maps with real multipliers},
  journal = {Transactions of the American Mathematical Society},
  volume  = {363},
  number  = {12},
  pages   = {6453--6463},
  year    = {2011},
  doi     = {10.1090/S0002-9947-2011-05308-0}
}

@article{Fatou1920,
  author  = {Fatou, Pierre},
  title   = {Sur les \'{e}quations fonctionnelles.
             {T}roisi\`eme m\'{e}moire},
  journal = {Bulletin de la Soci\'{e}t\'{e} Math\'{e}matique de France},
  volume  = {48},
  pages   = {208--314},
  year    = {1920},
  doi     = {10.24033/bsmf.1008}
}

@article{Inou,
  author  = {Inou, Hiroyuki},
  title   = {Extending local analytic conjugacies},
  journal = {Transactions of the American Mathematical Society},
  volume  = {363},
  number  = {1},
  pages   = {331--343},
  year    = {2011},
  doi     = {10.1090/S0002-9947-2010-05049-4}
}

@article{Muller,
  author  = {M\"uller, Peter},
  title   = {Decompositions of rational functions over real and complex
             numbers and a question about invariant curves},
  journal = {Illinois Journal of Mathematics},
  volume  = {59},
  number  = {4},
  pages   = {825--838},
  year    = {2015},
  doi     = {10.1215/ijm/1488186011}
}

@article{Ritt,
  author  = {Ritt, Joseph F.},
  title   = {Prime and composite polynomials},
  journal = {Transactions of the American Mathematical Society},
  volume  = {23},
  number  = {1},
  pages   = {51--66},
  year    = {1922},
  doi     = {10.1090/S0002-9947-1922-1501189-9}
}

@misc{MilnorIntro,
  author        = {Milnor, John},
  title         = {Dynamics in one complex variable: Introductory lectures},
  howpublished  = {Stony Brook IMS Preprint 1990/5},
  year          = {1991},
  note          = {Revised September 1991; arXiv:math/9201272},
  eprint        = {math/9201272},
  archivePrefix = {arXiv},
  primaryClass  = {math.DS},
  doi           = {10.48550/arXiv.math/9201272}
}

@article{Bowen1971,
  author  = {Bowen, Rufus},
  title   = {Entropy for group endomorphisms and homogeneous spaces},
  journal = {Transactions of the American Mathematical Society},
  volume  = {153},
  pages   = {401--414},
  year    = {1971},
  doi     = {10.1090/S0002-9947-1971-0274707-X}
}

@article{PakovichPolynomial,
  author  = {Pakovich, Fedor},
  title   = {Polynomial semiconjugacies, decompositions of iterations,
             and invariant curves},
  journal = {Annali della Scuola Normale Superiore di Pisa.
             Classe di Scienze (5)},
  volume  = {17},
  number  = {4},
  pages   = {1417--1446},
  year    = {2017},
  doi     = {10.2422/2036-2145.201510_019}
}

@incollection{MilnorLattes,
  author    = {Milnor, John},
  title     = {On {Latt\`es} maps},
  editor    = {Hjorth, Poul G. and Petersen, Carsten Lunde},
  booktitle = {Dynamics on the {Riemann} Sphere:
               A {Bodil Branner} Festschrift},
  publisher = {European Mathematical Society},
  address   = {Z\"urich},
  pages     = {9--43},
  year      = {2006},
  doi       = {10.4171/011-1/1}
}

\end{document}